\documentclass[12pt,twoside]{article}
\usepackage{amsfonts}
\usepackage{amssymb,amsmath}
\begin{document}
\begin{center} \textbf{The number of equations c = a+b satisfying the abc - conjecture} \vspace{12pt}
\\Constantin M. Petridi
\\ cpetridi@hotmail.com
\end{center}
\par
\vspace{2pt}
\begin{center}
\small{
\begin{tabular}{p{11cm}}
\textbf{Abstract} We prove that for a positive integer $c$ and
any given $\varepsilon$, $0<\varepsilon<1$, the number $N(c)$ of
equations $c=a+b$, $a<b$, with positive coprime integers $a$ and
$b$, which satisfy the inequality
$$c\,<\,R(c)^{\frac{\varepsilon}{1+\varepsilon}}R(a)^{\frac{1}{1+\varepsilon}}R(b)^{\frac{1}{1+\varepsilon}},$$
where $R(n)$ is the radical of $n$, is for $c\rightarrow\infty$
$$N(c)=(1-\varepsilon)\frac{\phi(c)}{2}+O\Bigl(\frac{\phi(c)}{2}\Bigr).$$
An analogue for the abc-conjecture inequality
$c<R(abc)^{1+\varepsilon}$ (without a constant factor) will also
be proved.
\end{tabular}
}
\end{center}
\par
\vspace{5pt} \textbf{1. Introduction}
\par\vspace{10pt}
 In our paper
arXiv:math/0511224v3[math.NT] 1 Mar 2006, we proved that for
positive coprime integers $a_{i},b_{i},c$, $1\leq i \leq
\frac{\varphi(c)}{2}$, satisfying $c=a_{i}+b_{i}$, $a_{i}<b_{i}$,
and for any given $\varepsilon>0$, there is a positive constant
$\kappa_{\varepsilon}$, effectively computable, depending on
$\varepsilon$, such that
$$\hspace{3cm}
\kappa_{\varepsilon}\,R(c)^{1-\varepsilon}\,c^{2}\;<\;\Bigl[
\prod_{1\leq i \leq
\frac{\varphi(c)}{2}}R(a_{i}b_{i}c)\Bigr]^{\frac{2}{\varphi(c)}}.\hspace{3cm}(1)
$$
Here $R(n)$ is the radical of $n$ and $\phi(n)$ is the Euler
totient function.
\par\vspace{10pt}
We shall use this result to estimate for a positive integer $c$
and any given $\varepsilon$, $0<\varepsilon<1$, the number of
equations $c=a+b$, $a\,<\,b$, with positive coprime integers $a$
and $b$, which satisfy the inequality
$$c\,<\,R(c)^{\frac{\varepsilon}{1+\varepsilon}}R(a)^{\frac{1}{1+\varepsilon}}R(b)^{\frac{1}{1+\varepsilon}}.$$
The analogous estimate for the abc-conjecture inequality
$$c\,<\,R(abc)^{1+\varepsilon},$$
follows as a consequence.
\par
\vspace{10pt} \textbf{2. Main Theorem}
\par\vspace{10pt}
\textbf{Theorem 1.} For a positive integer $c$ and any given
$\varepsilon$, $0<\varepsilon<1$, let $N(c)$, $1\leq N(c) \leq
\frac{\phi(c)}{2}$, be the number of equations $c =
a+b,\;a\,<\,b$ with coprime integers $a$ and $b$, which satisfy
the inequality
$$c<R(c)^{\frac{\varepsilon}{1+\varepsilon}}R(a)^{\frac{1}{1+\varepsilon}}R(b)^{\frac{1}{1+\varepsilon}}.$$
Then for $c\rightarrow\infty$
$$N(c)=(1-\varepsilon)\frac{\varphi(c)}{2}+O\Bigl(\frac{\varphi(c)}{2}\Bigr).$$
\par\vspace{10pt}
\textbf{Proof.} $N(c)$, has been defined as the number of
equations $c=a+b$, $a\,<\,b$ with positive coprime integers $a$
and $b$, satisfying
$$c\,<\,R(c)^{\frac{\varepsilon}{1+\varepsilon}}R(a)^{\frac{1}{1+\varepsilon}}R(b)^{\frac{1}{1+\varepsilon}},$$
which can also be written as
$$\hspace{4.3cm}R(c)^{1-\varepsilon}c^{1+\varepsilon}\,<\,R(cab).\hspace{4.3cm}(2)$$
On the other hand, because of
$c=a_{i}+b_{i},\;a_{i}\,<\,b_{i},\;(a_{i},b_{i})=1,\;1\leq i \leq
\frac{\phi(c)}{2}$, and $R(c)\leq c$, we have,
$$\hspace{3cm}R(a_{i}b_{i}c)=R(a_{i})R(b_{i})R(c)\,<\,R(c)c^{2}.\hspace{3cm}(3)$$
In the product $\Bigl[\prod_{1\leq i \leq
\frac{\varphi(c)}{2}}R(a_{i}b_{i}c)\Bigr]^{\frac{2}{\varphi(c)}}$,
therefore, because of (2), there are $N(c)$ factors, in some
order, which are greater than
$R(c)^{1-\varepsilon}c^{1+\varepsilon}$, but smaller than
$R(c)c^{2}$, as per (3). The remaining $\frac{\phi(c)}{2}-N(c)$
factors, according to same definition of $N(c)$, are all smaller
than $R(c)^{1-\varepsilon}c^{1+\varepsilon}$.
\par\vspace{10pt}
In view of this and of (1), we deduce that
$$\kappa_{\varepsilon}R(c)^{1-\varepsilon}c^{2}\,<\,
\Bigl[(R(c)c^{2})^{N(c)}\Bigr]^{\frac{2}{\phi(c)}}\Bigl[(R(c)^{1-\varepsilon}c^{1+\varepsilon})^{\frac{\phi(c)}{2}-N(c)}\Bigr]^{\frac{2}{\phi(c)}}.$$
Simplifying, we get
$$\kappa_{\varepsilon}R(c)^{1-\varepsilon}c^{2}\,<\,\Bigl(R(c)c^{2}\Bigl)^{\frac{2}{\phi(c)}N(c)}\Bigl(R(c)^{1-\varepsilon}c^{1+\varepsilon}\Bigr)\Bigl(R(c)^{\varepsilon -1}c^{-1-\varepsilon}\Bigr)^{\frac{2}{\phi(c)}N(c)},$$
$$\kappa_{\varepsilon}c^{1-\varepsilon}\,<\Bigl(R(c)^{\varepsilon}\,c^{1-\varepsilon}\Bigr)^{\frac{2}{\phi(c)}N(c)}.$$
We now take the logarithms of both sides to obtain
$$\log\kappa_{\varepsilon}+(1-\varepsilon)\log c<\Bigl(\varepsilon\log R(c)+(1-\varepsilon)\log c\Bigr)\frac{2}{\phi(c)}N(c).$$
Dividing by $(\varepsilon\log R(c)+(1-\varepsilon)\log c)\,>\,0$
and noting that $\frac{2}{\phi(c)}N(c)\leq 1$, we get
$$\frac{\log\kappa_{\varepsilon}+(1-\varepsilon)\log c}{\varepsilon
\log R(c)+(1-\varepsilon)\log c}\,<\,\frac{2}{\phi(c)}N(c)\leq
1.$$
\par\vspace{10pt}
Since $\log R(c)$ is less than $\log c$, we conclude that
$$\frac{\log\kappa_{\varepsilon}+(1-\varepsilon)\log c}{\log
c}<\frac{2}{\phi(c)}N(c)\,\leq\,1.$$ Thus
$$\frac{\log\kappa_{\varepsilon}}{\log
c}+(1-\varepsilon)\,<\,\frac{2}{\phi(c)}N(c)\,\leq\,1,$$ or,
written otherwise,
$$\frac{\log\kappa_{\varepsilon}}{\log
c}\,<\,\frac{2}{\phi(c)}N(c)-(1-\varepsilon)\,\leq\,\varepsilon.$$
By letting $c\rightarrow\infty$, this gives
$$N(c)=(1-\varepsilon)\frac{\phi(c)}{2}+O\Bigl(\frac{\phi(c)}{2}\Bigr),$$
as claimed by Theorem 1.
\par
\vspace{20pt} \textbf{3. Analogue for the abc-conjecture}
\par\vspace{10pt}
\textbf{Theorem 2.} For a positive integer $c$ and any given
$\varepsilon$, $0<\varepsilon<1$, let $N_{1}(c)$, $1\leq N_{1}(c)
\leq \frac{\phi(c)}{2}$, be the number of equations $\,c =
a+b,\;a\,<\,b$ with coprime integers $a$ and $b$, which satisfy
the inequality
$$c\,<\,R(c)^{1+\varepsilon}R(a)^{1+\varepsilon}R(b)^{1+\varepsilon}.$$
Then for $c\rightarrow\infty$
$$N_{1}(c)=(1-\varepsilon)\frac{\varphi(c)}{2}+O\Bigl(\frac{\varphi(c)}{2}\Bigr).$$
\par
\textbf{Proof.} Since
$1+\varepsilon\,>\,\frac{\varepsilon}{1+\varepsilon}$ and
$1+\varepsilon\,>\,\frac{1}{1+\varepsilon}$, we have
$$c\,<\,R(c)^{\frac{\varepsilon}{1+\varepsilon}}R(a)^{\frac{1}{1+\varepsilon}}R(b)^{\frac{1}{1+\varepsilon}}\,<\,R(c)^{1+\varepsilon}R(a)^{1+\varepsilon}R(b)^{1+\varepsilon}.$$
This means that the set of equations $\,c = a+b,\;a\,<\,b$ with
coprime integers $a$ and $b$, satisfying Theorem 1, does, a
fortiory, also satisfy Theorem 2.
\par\vspace{10pt}
As a consequence $N_{1}(c)\geq N(c)$, and as
$N(c)=(1-\varepsilon)\frac{\varphi(c)}{2}+O\Bigl(\frac{\varphi(c)}{2}\Bigr)$,
according to Theorem 1, it also follows that
$$N_{1}(c)=(1-\varepsilon)\frac{\varphi(c)}{2}+O\Bigl(\frac{\varphi(c)}{2}\Bigr),$$
which proves the Theorem 2.
\par\vspace{20pt}
In a next paper we examine for which functions $H(x,y,z)$, the
inequality
$$c\,<\,H\bigl(R(c),R(a),R(b)\bigr),$$
in combination with
$$
\kappa_{\varepsilon}\,R(c)^{1-\varepsilon}\,c^{2}\;<\;\Bigl[
\prod_{1\leq i \leq
\frac{\varphi(c)}{2}}R(a_{i}b_{i}c)\Bigr]^{\frac{2}{\varphi(c)}},
$$
can yield substantial results.
\par
 \vspace{20pt} \textbf{Acknowledgment.} I am indebted
to Peter Krikelis, Department of Mathematics, University of
Athens, for his unfailing assistance.

\end{document}